\documentclass[letterpaper, 10 pt, conference]{ieeeconf}  

\IEEEoverridecommandlockouts                              

\usepackage{graphics} 
\usepackage{epsfig} 
\usepackage{mathptmx} 
\usepackage{times} 
\usepackage{amsmath} 
\usepackage{amssymb}  
\usepackage{tabularx}
\usepackage{algorithm}
\usepackage{algorithmic}
\usepackage{subcaption}
\usepackage{url}

\DeclareFontEncoding{LS1}{}{}
\DeclareFontSubstitution{LS1}{stix}{m}{n}
\DeclareSymbolFont{stixletters}{LS1}{stix}{m}{it}
\DeclareMathAccent{\cev}{\mathord}{stixletters}{"91}

\title{\LARGE \bf
Time-reversal solution of BSDEs in stochastic optimal control:\\ a linear quadratic study
}

\author{Yuhang Mei$^{1}$ and Amirhossein Taghvaei$^{1}$ 
\thanks{$^{1}$Y. Mei and A. Taghvaei are with the William E. Boeing department of Aeronautics and Astronautics at University of Washington, Seattle;
        {\tt\small yuhangm@uw.edu, amirtag@uw.edu}. This research is supported by the National Science Foundation (NSF) award EPCN-2347358. }%
 }

\newcommand{\ud}{\mathrm{d}}
\def\Re{\mathbb{R}}
\def\R{\mathbb{R}}
\def\argmin{\mathop{\text{\rm arg\,min}}}

\newcommand{\tr}{\mbox{tr}}

\newtheorem{remark}{Remark}

\newtheorem{assumption}{Assumption}

\newcommand{\trace}{\text{Tr}}

\newcounter{rmnum}

\newcounter{anum}

\begin{document}

\maketitle
\thispagestyle{empty}
\pagestyle{empty}

\begin{abstract}
This paper addresses the numerical solution of backward stochastic differential equations (BSDEs) arising in stochastic optimal control. Specifically, we investigate two BSDEs: one derived from the Hamilton-Jacobi-Bellman equation and the other from the stochastic maximum principle. For both formulations, we analyze and compare two numerical methods. The first utilizes the least-squares Monte-Carlo (LSMC) approach for approximating conditional expectations, while the second leverages a time-reversal (TR) of diffusion processes. Although both methods extend to nonlinear settings, our focus is on the linear-quadratic case, where analytical solutions provide a benchmark. Numerical results demonstrate the superior accuracy and efficiency of the TR approach across both BSDE representations, highlighting its potential for broader applications in stochastic control.   
\end{abstract}

\section{INTRODUCTION}
	Stochastic optimal control (SOC) is central to applications from various disciplines, including   finance~\cite{pham2009continuous}, stochastic thermodynamics~\cite{sekimoto2010stochastic,peliti2021stochastic}, robotic navigation~\cite{shah2012stochastic,exarchos2019optimal}, and diffusion models~\cite{ramesh2022hierarchical,rombach2022high}. BSDE solutions of stochastic optimal control problem have gained attention in the recent years ~\cite{exarchos2016learning,exarchos2018stochastic, pereira2019learning,han2018solving,archibald2020stochastic}
 due to the fundamental challenges of directly solving the Hamilton-Jacobi-Bellman (HJB) partial differential equation (PDE).   Namely, there are two BSDEs that arise in the context of stochastic optimal control. The first one is based on the BSDE formulation of the HJB equation,  derived from the application of the nonlinear Feynman-Kac lemma~\cite{exarchos2016learning,exarchos2018stochastic}. And the second one is obtained from the application of stochastic maximum principle~\cite{yong1999stochastic}. In the first formulation, the backward process identifies the value function evaluated along the trajectory of the state process. In the second formulation, the backward process is the co-state, which is equal to to the derivative of the value function, if the derivative exists,  evaluated along the state trajectory. 
 In order to differentiate the BSDEs, the first one is referred to as the value function BSDE, while the second one  is referred to as the co-state BSDE. 
 
 In both settings, numerical solution of the BSDE is faced with the challenge of finding a forward adapted process  that satisfies a terminal constraint. 
  There are several PDE-based approaches~\cite{peng1991probabilistic,ma1994solving,ma2002numerical} that solve the BSDE, but they also suffer from  limitations  that  are similar to solving the HJB equation. 

  In this paper, we study two sampling-based approaches that numerically approximate the solution to the BSDE. The first approach is based 
   on expressing the backward process as conditional expectation with respect to the forward filtration~\cite{zhang2004numerical,bouchard2004discrete,zhao2006new,bender2007forward,zhang2013sparse}. In particular,  we focus on~\cite{exarchos2016learning,exarchos2018stochastic} where the least-squares Monte-Carlo (LSMC) approach is proposed for numerically approximating the conditional expectations. 
   The second approach is based on time-reversal (TR) of diffusion~\cite{anderson1982reverse,haussmann1986time,follmer2005entropy,cattiaux2021time}. Specifically, we study the TR approach, presented as part of our prior work~\cite{taghvaei2024time},  designed to numerically solve the co-state BSDE. The TR formulation reverses the direction of the filteration, thus facilitating the numerical simulation of the backward process, by simply initiating it at terminal time and  integrating it backward with the Euler-Maruyama scheme. The TR comes with the additional cost of approximating the so-called F\"ollmer's drift or score function, for which efficient numerical methods have been developed in the context of diffusion models~\cite{song2019generative,song2020score,ho2020denoising}. 

The objective of the paper is to study and compare the LSMC and TR approach for solving the value function and co-state BSDEs. In order to do so, we present both algorithms within a unified framework that facilitates the comparison. In doing so, we provide an elementary derivation of the TR formulation for BSDEs, that is more general than the one presented in the prior work, and applies to value function BSDE. We provide extensive numerical results in the linear quadratic (LQ) setting, analyzing the effect of time-discretization, sample size, and problem dimension. 
Our results demonstrate that solving the co-state BSDE is generally more numerically stable compared to the value function BSDE, and the TR approach gives significantly better accuracy compared to the LSMC. These findings highlight the potential of the TR approach and motivate future research, particularly focusing on  the nonlinear setting.

The outline of the paper is as follows: Section~\ref{sec:problem} presents the stochastic optimal control problem and the value function and co-state BSDEs. Section~\ref{sec:algorithms} contains the two numerical algorithms, LSMC and TR, that are studied in this paper. Finally, the numerical results and concluding remarks are presented in Section~\ref{sec:numerics} and~\ref{sec:numerics}, respectively.

\section{Problem formulation and background}\label{sec:problem}
In this paper, we are interested in control systems that are governed by the stochastic differential equation (SDE)
\begin{equation}\label{eq:dynamic-SDE}
    \ud X_t = a(X_t,U_t) \ud t + \sigma(X_t) \ud W_t, \quad X_0 \sim p_0 ,
\end{equation}
where $X:=(X_t)_{t\geq 0}$  is the $n$-dimensional state trajectory, $ U:=( U_t)_{t\geq 0}$  is the $m$-dimensional control input,  $ W:= (W_t)_{t\geq 0}$  is the standard $n$-dimensional Wiener process, and $p_0$  is the probability distribution of the initial condition $X_0$. The functions $a: \R^n \times \R^m \rightarrow \R^n$ and $\sigma: \R^n \rightarrow \R^{n\times n}$ are assumed to be smooth and globally Lipschitz. The control input $U$ is assumed to be adapted to the filtration $\mathcal{F}_t=\sigma(W_s;\,0\leq s\leq t)$ generated by the Wiener process $W$. This constraint is denoted by $U_t \in \mathcal F_t$ and highlights the non-anticipative nature of the control. The cost function, associated with the control system~\eqref{eq:dynamic-SDE}, is
\begin{equation}\label{eq:cost}
    J(t, x; U) := \mathbb{E} \left [ \int_t^T \ell(X_s,U_s) \ud s + \ell_f(X_T) \,\bigg| \,X_{t}=x\right],
\end{equation}
for all $t \in [0,T]$, $x \in \Re^n$, and $U_t \in \mathcal F_t$, 
where the running cost $\ell:\R^n \times \R^m \rightarrow \R$, and the terminal cost $\ell_f: \R^n \rightarrow \R$ are assumed to be smooth, and the horizon $T\in(0,\infty)$ is assumed to be finite. The stochastic optimal 
problem is to solve the optimization problem
\begin{equation}\label{eq:SOC}
    \inf_{U_t \in \mathcal F_t}\, \mathbb E[J(0,X_0;U)], ~~s.t. ~ \eqref{eq:dynamic-SDE}.
\end{equation}

In the following two subsections, we present two BSDEs that characterize the solution to the SOC problem~\eqref{eq:SOC}. 



\subsection{Value function BSDE}
The value function is defined according to
\begin{equation}
    V(t, x) = \inf_{U_t \in \mathcal F_t} J(t, x; U)
    ,\quad \forall ~(t,x) \in  [0,T] \times \Re^n. 
\end{equation}
The dynamic programming principle implies that the value function solves the HJB PDE~\cite[Prop. 3.5]{yong1999stochastic}
\begin{align}
    \nonumber
   \partial_t &V (t,x) + \inf_u \{ \partial_x V(t,x)^\top a(x,u) + \ell(x,u) \} \\ &+\frac{1}{2}\trace( \partial_{xx} V(t,x)\sigma(x) \sigma(x)^\top) = 0,\quad 
    V(T,x) = \ell_f(x),\label{eq:HJB}
\end{align} 
for all $x\in \Re^n$ and $t\in[0,T]$. 
Given the value function, 
the optimal control input is 
\begin{equation}\label{eq:optimal-control-value}
 U_t 
 =\argmin_u \, \{ \partial_x V(t,X_t)^\top a(X_t,u) + \ell(X_t,u) \} 
\end{equation}
for any realization $X_t$ of the state process (it is  assumed that the minimization~\eqref{eq:optimal-control-value} has a unique solution). 

The HJB equation  leads to a BSDE under the following assumption about the model and cost.
\begin{assumption}\label{assum:control-affine}
    The dynamic model admits a control affine form $a(x,u) = \tilde a(x) + \sigma(x)\tilde B(x)u$ and the cost function  is separable and quadratic in control, i.e. $\ell(x,u) = \tilde \ell(x) + \frac{1}{2}u^\top R u$ for a positive definite matrix $R$. 
\end{assumption}

The BSDE  is formally constructed as follows. 
Let $X_t$ solve~\eqref{eq:dynamic-SDE} with some  control input $U_t \in \mathcal F_t$. Application of the It\^o rule to $V(t,X_t)$, along with the HJB equation~\eqref{eq:HJB}, concludes that 
\begin{align*}
    \ud &V(t,X_t) = \partial_x V(t,X_t)^\top a(X_t,U_t) \ud t + \partial_x V(t,X_t)^\top \sigma(X_t)\ud W_t  \\-& \inf_u \{ \partial_x V(t,X_t)^\top a(X_t,u) + \ell(X_t,u) \}\ud t,~ V(T,X_T)=\ell_f (X_T).
\end{align*}
Using Assumption~\ref{assum:control-affine}, and 
defining the pair
\begin{equation}\label{eq:Yval-Zval}
    (Y_t,Z_t)\mathrel{:=}(V(t,X_t), \sigma(X_t)^\top \partial_x V(t,X_t))
\end{equation} concludes the BSDE 
\begin{align}\label{eq:BSDE-value-function}
    -\ud  Y_t = h(X_t,U_t,Y_t,Z_t)\ud t - Z_t^\top \ud W_t,\quad Y_T = \ell_f (X_T),
\end{align}
where 
\begin{align}\label{eq:h}
    h(x,u,y,z):= \tilde \ell(x) - \frac{1}{2} z^\top \tilde B(x) R^{-1}\tilde B(x)^\top z- z^\top \tilde B(x) u,  
\end{align}
for $(x,u,y,z)\in \Re^n \times \Re^m \times \Re \times \Re^n$. The optimal control input is identified in terms of the solution with the formula
\begin{equation}\label{eq:optimal-control-value-affine}
    U^*_t = -R^{-1} \tilde B(X_t)^\top Z_t. 
\end{equation}
A rigorous analysis of the BSDE~\eqref{eq:BSDE-value-function} and its relation to the nonlinear Feynman-Kac lemma appears in~\cite{exarchos2018stochastic}. 
\subsection{Co-state BSDE}
The Hamiltonian function for the SOC problem~\eqref{eq:SOC} is defined according to 
\begin{equation}
    H(x,u,y,z) := \ell (x,u) + y^\top a(x,u) + \trace(z \sigma(x)), 
\end{equation}
for $(x,u,y,z) \in \Re^n\times \Re^m\times \Re^n\times \Re^{n\times n}$.
The Hamiltonian function is used to introduce 
the BSDE 
\begin{align}
    -\ud Y_t &= \partial_x H (X_t, U_t, Y_t, Z_t) \ud t - Z_t^\top \ud W_t, ~ Y_T = \partial_x \ell_f(X_T),\label{eq:MP-BSDE}
\end{align}
for any control input $U_t \in \mathcal F_t$, where the pair $(Y,Z):=(Y_t,Z_t)_{t\geq 0}$ are known as adjoint processes. 
According to the stochastic maximum principle, the optimal control input $U_t$ satisfies
\begin{equation}\label{eq:min-U}
    U_t =\argmin_u\,H(X_t, u,Y_t, Z_t),\quad \forall t \in[0,T].
\end{equation}
where $(X,Y,Z)$ are the solutions to~\eqref{eq:dynamic-SDE}-\eqref{eq:MP-BSDE}~\cite[Thm. 3.2, Case 1]{yong1999stochastic}. 

The solution to the BSDE \eqref{eq:MP-BSDE} may be expressed in terms of the value function, when it is differentiable, according to 
\begin{equation}
    \label{eq:Y-Z-MP}
    (Y_t,Z_t)=(\partial_x V(t,X_t),\sigma(X_t)^\top \partial_{xx} V(t,X_t)).
\end{equation}


\begin{remark}
The stochastic maximum principle does not require   the control-affine model and separable cost Assumption~\ref{assum:control-affine}. 
\end{remark}


\subsection{Linear quadratic case}

A special class of stochastic optimal control problems is the linear quadratic (LQ) case where the dynamic model and cost functions take the special form
\begin{align*}
    &a(x,u)= A x + B u,\quad \sigma(x) = \sigma,\\
    &\ell(x, u) = \frac{1}{2} x^\top Q x + \frac{1}{2} u^\top R u, \quad \ell_f (x) = \frac{1}{2} x^\top Q_f x,
\end{align*}
where $A$, $B$, $\sigma$, $Q$, $R$, and $Q_f$ are matrices of appropriate dimensions. In the LQ case, it is further assumed that $Q$ and $Q_f$ are symmetric positive semidefinite matrices, while $R$ is a symmetric positive definite matrix. 

The LQ case is an important benchmark because the exact solution is explicitly known in this case. 
In particular, the value function $V(t,x)=\frac{1}{2}x^\top G_t x + g_t$ where $G_t$ and $g_t$  solve the differential equations 
\begin{align}\label{eq:Ricatti}
-\dot{G}_t &= G_tA + A^\top G_t + Q - G_tBR^{-1}B^\top G_t,\quad G_T = Q_f\\
-\dot g_t &= \frac{1}{2} \trace(\sigma \sigma^\top G_t) \quad g_T = 0 \nonumber 
\end{align}
The optimal control law is given by
\begin{equation}
    U_t=- R^{-1} B^\top G_t X_t. 
\end{equation}
The solution to the BSDEs~\eqref{eq:BSDE-value-function} and \eqref{eq:MP-BSDE} is also explicitly known and given by 
\begin{align*}
   \text{BSDE~\eqref{eq:BSDE-value-function}:}& \quad (Y_t,Z_t)= (\frac{1}{2}X_t^\top G_t X_t, \sigma^\top G_t X_t)\\
\text{BSDE~\eqref{eq:MP-BSDE}:}& \quad (Y_t,Z_t)= (G_t X_t, \sigma^\top G_t ). 
\end{align*}

\section{Algorithms}\label{sec:algorithms}
In this section, we present two algorithms  that aim to solve the BSDEs~\eqref{eq:BSDE-value-function} and \eqref{eq:MP-BSDE}. In both cases, we assume the control input is generated through  a feedback control law  $U_t=k(t,X_t)$ for some given function $k:[0,T]\times \Re^n \to \Re^m$.  Later in Section~\ref{sec:numerics}, we explain how the solution to the BSDEs is used to iteratively update the feedback control law until convergence to the optimal control law. 

In order to facilitate the presentation, the BSDEs ~\eqref{eq:BSDE-value-function} and \eqref{eq:MP-BSDE} are expressed in the unified form
\begin{equation}\label{eq:BSDE-general}
		-\ud Y_t = g(t,X_t,Y_t,Z_t)\ud t - Z_t^\top \ud W_t,\quad Y_T = g_f(X_T),
	\end{equation}
where, depending on the BSDE, $g$ and $g_f$ take the following forms:
\begin{align*}
	\text{BSDE~\eqref{eq:BSDE-value-function}:}&~g(t,x,y,z)=h(x,k(t,x),y,z),\quad g_f(x)=\ell_f(x), \\
	\text{BSDE~\eqref{eq:MP-BSDE}:}&~ 	
	g(t,x,y,z)=\partial_x H(x,k(t,x),y,z),~ g_f(x)=\partial_x \ell_f(x).
\end{align*}

Furthermore,  we assume the PDE
\begin{align}\nonumber
	&\partial_t \phi (t,x) + \partial_x \phi(t,x)^\top a(x,k(t,x))  + \frac{1}{2}\tr(\partial_{xx}\phi(t,x)D(x)) \\&+  g(t,x,\phi(t,x),\sigma(x)^\top\partial_x \phi(t,x)) = 0,~ \phi(T,x) = g_f(x),\label{eq:phi-pde}
\end{align}
admits a smooth solution $\phi$. Under this assumption, the solution to the BSDE~\eqref{eq:BSDE-general} is given by the formula 
\begin{align}\label{eq:Y-Z-phi}
	(Y_t,Z_t)=(\phi(t,X_t),\sigma(X_t)^\top\partial_x\phi(t,X_t)). 
\end{align} 
For the case of BSDE~\eqref{eq:BSDE-value-function}, $\phi(t,x)=V(t,x)$ solves the PDE~\eqref{eq:phi-pde}, while for BSDE~\eqref{eq:MP-BSDE}, $\phi(t,x)=\partial_x V(t,x)$ solves the PDE~\eqref{eq:phi-pde}, where $V(t,x)$ is the value function.

The algorithms that are presented below aim at approximating the function $\phi(t,x)$, not by solving the PDE~\eqref{eq:phi-pde}, but through a sampling-based approach. 

%
\subsection{Least-Squares-Monte-Carlo (LSMC)} 
This section presents the LSMC algorithm, introduced in~\cite{exarchos2016learning,exarchos2018stochastic}, for approximating the conditional expectation that identifies the solution to the BSDE~\eqref{eq:BSDE-general}. 

Integrating the BSDE over the interval $[s,t] \subset [0,T]$ and taking the conditional expectation with respect to $\mathcal F_s$ yields 
\begin{align*}
	Y_s = \mathbb E[Y_t + \int_s^t g(\tau,X_\tau,Y_\tau, Z_\tau) \ud \tau\vert \mathcal F_s]. 
\end{align*}
Then, the minimum mean-squared-error property of conditional expectation and the relationship  $Y_s=\phi(s,X_s)$ concludes the following optimization problem for $\phi(s,\cdot)$
\begin{align*}
	\phi(s,\cdot)= \argmin_{\phi} \mathbb E[\|Y_t + \int_s^t g(\tau,X_\tau,Y_\tau, Z_\tau) \ud \tau- \phi(X_s)\|^2]. 
\end{align*}
Finally, by discretizing the time domain with step-size $\Delta t>0$ and generating $N$ independent realizations of the state trajectory $\{X^i_t\}_{i=1}^N$, we end up with a recursive update formula for the function $\phi$, such that, given $\phi(t,\cdot)$, the function $\phi(t-\Delta t,\cdot)$ is obtained by solving the minimization problem 
  \begin{align*}
 	\phi(t-\Delta t,\cdot)= \argmin_{\phi \in \Phi}\, \sum_{i=1}^N \|Y^i_t+  g(t,X^i_t,Y^i_t, Z^i_t) \Delta t - \phi(X^i_{t-\Delta t}) \|^2,
 \end{align*}
where $(Y^i_t,Z^i_t)=(\phi(t,X^i_t),\sigma(X^i_t)^\top \partial_x \phi(t,X^i_t))$, and $\Phi$ is a given parametric class of functions. 
A detailed description of the method is presented in Algorithm~\ref{alg:LSMC}. 

\begin{algorithm}[t]
	\caption{Least-squares-Monte-Carlo (LSMC)} 
	\begin{algorithmic}[1]
		\STATE \textbf{Input:} sample size $N$, step-size $\Delta t$, control law $k(t,x)$, function class $\Phi$.   
		\STATE \ $\{X_0^i\}_{i=1}^N \sim p_0$
		\FOR{$t\in \{0,\Delta t,\ldots, T-\Delta t\}$}
		\STATE $U^i_t= k(t, X^i_t)$ and $\{\Delta W^i_t\}_{i=1}^N\sim N(0,\Delta t I_n)$
		\STATE $X^{i}_{t+\Delta t} \!=\! X^i_{t}+a(X^{i}_{t},U^i_{t})\Delta t  + \sigma(X^i_t) \Delta W^i_t$
		\ENDFOR
		\STATE $Y^i_T = g_f(X^i_T)$, $Z^i_T = \sigma(X^i_T)^\top\partial_x g_f(X^i_T)$, and $\phi(T,\cdot)=g_f(\cdot)$
		\FOR{$t \in\{T,T-\Delta t, \ldots,\Delta t\}$}
  \STATE $\Delta Y^i_t=g(X^i_t,U^i_t,Y^i_t,Z^i_t)\Delta t$
		\STATE $\phi(t-\Delta t,\cdot)= \argmin_{\phi \in \Phi}\,\sum_{i=1}^N \left\|Y^i_{t} \!+\! \Delta Y^i_t \!-\! \phi(X^i_{t-\Delta t})\right\|^2$
		\STATE $Y^i_{t-\Delta t} = \phi(t-\Delta t,X^i_{t-\Delta t})$ 
		\STATE $Z^i_{t-\Delta t} = \sigma(X^i_{t-\Delta t})^\top\partial_x \phi(t-\Delta t,X^i_{t-\Delta t})$
		\ENDFOR
		\STATE \textbf{Output:}  $\{\phi(t,\cdot)\}_{t\in \{0,\Delta t,\ldots,T\} }$
	\end{algorithmic}
	\label{alg:LSMC}
\end{algorithm}



 \subsection{Time-Reversal (TR)}
 This section presents a derivation of the TR approach for BSDEs, that is more general and elementary compared to the one introduced in the prior work~\cite{taghvaei2024time}. In particular,~\cite{taghvaei2024time} uses the stochastic maximum principle for Mckean-Vlasov SDEs to derive the time-reversal formulation of only~\eqref{eq:MP-BSDE}, while, here we use \eqref{eq:phi-pde} and It\^o's formula to derive  the time-reversal formulation for both \eqref{eq:BSDE-value-function} and \eqref{eq:MP-BSDE}.
 
 The SDE~\eqref{eq:dynamic-SDE}, with feedback control $U_t=k(t,X_t)$,  
 admits the time-reversal formulation
 \begin{align} \label{eq:dyn-sde-reverse}
 	\ud \cev X_t\! &=\! a (\cev X_t,k(t,\cev X_t)) \ud t \!+\! b(t,\cev X_t) \ud t \!+\! \sigma(\cev X_t) \cev \ud  W_t, ~ \cev X_T\! \sim\! p_{T}
 \end{align}
 where $\cev  \ud W_t$ denotes backward stochastic integration~\cite[Sec. 2]{pardoux1987two}, $p_{t}(x)$ is the probability density function of $X_t$,
 \[b(t,x) :=- \frac{1}{p_{t}(x)}\partial_x (D(x)^\top p_{t}(x))\] 
 is known as the {F\"ollmer's drift}  term or score function, and $D(x) := \sigma(x)\sigma(x)^\top$. The solution $\cev X_t$ to the time-reversed SDE~\eqref{eq:dyn-sde-reverse} is adapted to the backward filteration $\cev {\mathcal F}_t:=\sigma(W_T-W_s;\,t\leq s\leq T)$. Moreover, the solution is equal, in probability, to the solution to the SDE~\eqref{eq:dynamic-SDE}, i.e. the probability distribution of the path $\{\cev X_t;\,0\leq t\leq T\}$ is equal to the  probability distribution of the path $\{X_t;\,0\leq t\leq T\}$~\cite{cattiaux2021time}.

Applying the integration by parts and following \cite[Thm. 1]{hyvarinen2005estimation}, the score function $b(t,x)$ is the solution to the minimization
 \begin{align*}
 	\min_b\,\mathbb E\left[\frac{1}{2}\|b(X_t)\|^2 - \trace(D(X_t)\partial_x b(X_t)) \right].
 \end{align*}
This result is used for the numerical  approximation of the score function. In particular,  given $N$ realizations of the state $\{X^i_t\}_{i=1}^N$, the score function is approximated by solving the optimization problem
\begin{align}\label{eq:b-min}
    \min_{b\in \Psi}\,\frac{1}{N} \sum_{i=1}^N \left[\frac{1}{2}\|b(X^i_t)\|^2 + \trace(D(X^i_t)\partial_x b(X^i_t)) \right],
\end{align}
where $\Psi$ is a given class of functions. 
 
  \begin{figure*}[t]
	\centering
	\includegraphics[width=0.22\hsize,trim={8pt 0pt 5pt 0 pt},clip]{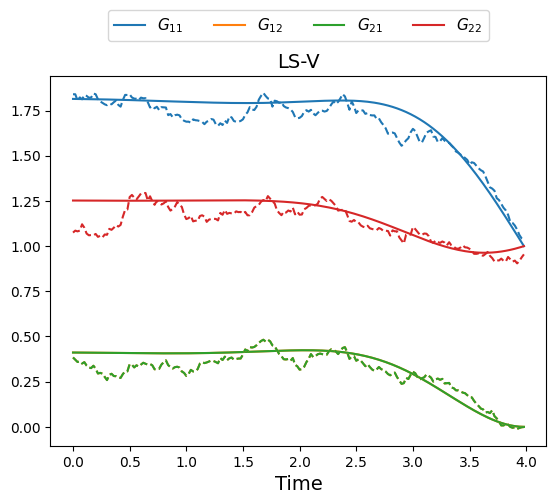} \hfill 
	\includegraphics[width=0.22\hsize,trim={8pt 0pt 5pt 0 pt},clip]{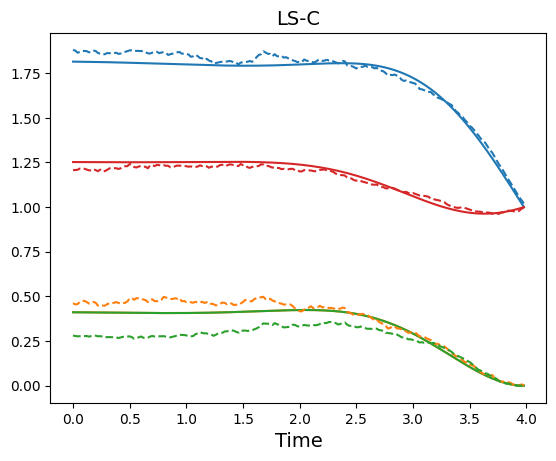} \hfill 
	\includegraphics[width=0.22\hsize,trim={8pt 0pt 5pt 0 pt},clip]{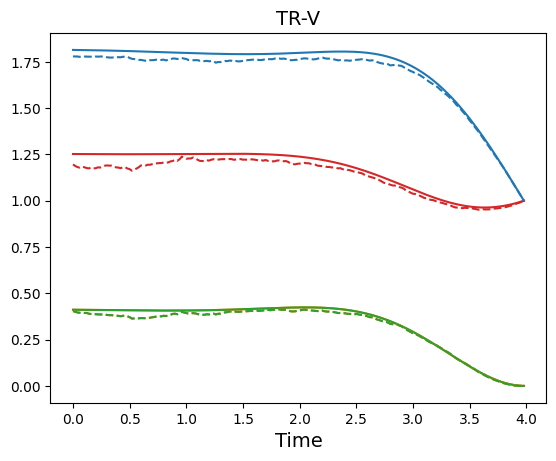}\hfill
        \includegraphics[width=0.22\hsize,trim={8pt 0pt 5pt 0 pt},clip]{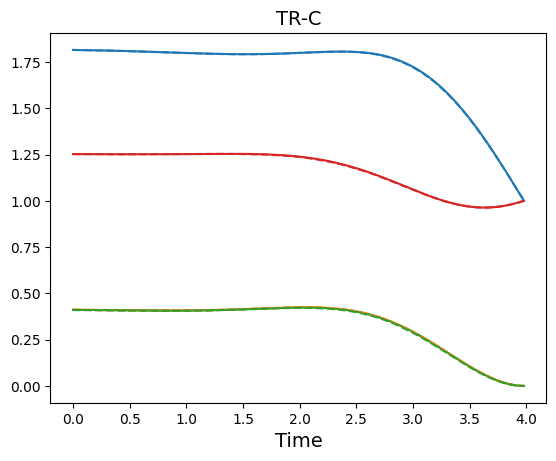}\\
	\hspace{10pt}(a) \hspace{115pt} (b)  \hspace{115pt}  (c) \hspace{110pt} (d)
	\caption{Numerical result for Section~\ref{sec:numerics-accuracy}: The entries of the matrix $G_t$, obtained from the four algorithms (a) LSMC-V, (b) LSMC-C, (c) TR-V, and (d) TR-C, in comparison to their exact values. The solid line denotes the exact solution and the dotted line denotes the numerical approximation from algorithms. The  time horizon $T=4$, time step-size $\Delta t = 0.02$, and sample size $N=2000$.}
	\label{fig:g_matrices}
\end{figure*}
The time-reversed process $\cev X_t$ is used to construct the time-reversal of the BSDE~\eqref{eq:BSDE-general} as follows. Application of the 
 It\^o rule to $\phi(t,\cev X_t)$, where $\phi$ is the solution to the PDE~\eqref{eq:phi-pde}, concludes 
 \begin{align*}
 	\ud \phi(t,\cev X_t) = -&g(t,\cev X_t,\phi(t,\cev X_t),\sigma(\cev X_t)^\top\partial_{x} \phi(t,\cev X_t)) \ud t \\
 	-& c(t,\cev X_t)\ud t+ \partial_x \phi(t,\cev X_t)^\top\sigma(\cev X_t) \cev \ud W_t,
 \end{align*}
where $c(t,x):=\trace(D(x)\partial_{xx} \phi(t,x)) -\partial_x \phi(t,x)^\top b(t,x)$. 
 Defining $(\cev Y_t,\cev Z_t)=(\phi(t,\cev X_t),\sigma(\cev X_t)^\top \partial_{x}\phi(t,\cev X_t))$ concludes the relationship
 \begin{equation}\label{eq:BSDE-general-reversed}
 	-\ud \cev Y_t = g(t,\cev X_t,\cev Y_t,\cev Z_t)\ud t - \cev Z_t^\top \cev \ud W_t + c(t,\cev X_t)\ud t,\quad \cev Y_T = g_f(\cev X_T)
 \end{equation}
 which is viewed as the time-reversal of the BSDE~\eqref{eq:BSDE-general}.  By construction, the solution $(\cev Y_t,\cev Z_t)$ is adapted to the backward filteration $\cev {\mathcal F}_t$, and the process $(Y,Z)=\{(Y_t,Z_t);\,0\leq t\leq T\}$ has the same probability law as  $(\cev Y,\cev Z)=\{(\cev Y_t,\cev Z_t);\,0\leq t\leq T\}$. Moreover, the relationship  $\cev Y_t=\phi(t,\cev X_t)$ implies the following minimization problem to numerically approximate $\phi(t,x)$: 
 \begin{equation}\label{eq:phi-min}
 	\min_{\phi \in \Phi}\, \sum_{i=1}^N \| \cev Y^i_t- \phi(t,\cev X^i_t)\|
 \end{equation}
 where $\{(\cev X^i_t,\cev Y^i_t)\}_{i=1}^N$ are realizations of~\eqref{eq:dyn-sde-reverse}-\eqref{eq:BSDE-general-reversed}. 

The time-reversed BSDE~\eqref{eq:BSDE-general-reversed} is used to develop a numerical algorithm to solve BSDEs. Similar to LSMC, the algorithm  starts with simulating $N$ independent realizations of the state process $\{X^i_t\}_{i=1}^N$ over a discretization of the time-domain. However,  the algorithm involves the additional step of approximating the score function by solving the minimization~\eqref{eq:b-min}. Then, the time-reversed state and adjoint processes~\eqref{eq:dyn-sde-reverse}-\eqref{eq:BSDE-general-reversed} are simulated backward, while simultaneously  the minimization problem~\eqref{eq:phi-min} is solved to to obtain  the function $\phi$. The details are presented in Algorithm~\ref{alg:TR}.

\begin{algorithm}[h]
\caption{Time-reversal (TR)} 
\begin{algorithmic}[1]
\STATE \textbf{Input:}  sample size $N$, step-size $\Delta t$, control law $k(t,x)$, function class $\Phi$ and $\Psi$.   
\STATE  $\{X_0^i\}_{i=1}^N \sim p_0$
\FOR{$t\in \{0,\Delta t,\ldots, T-\Delta t\}$}
\STATE $U^i_t= k(t, X^i_t)$ and $\{\Delta W^i_t\}_{i=1}^N\sim N(0,\Delta t I_n)$
		\STATE $X^{i}_{t+\Delta t} \!=\! X^i_{t}+a(X^{i}_{t},U^i_{t})\Delta t  + \sigma(X^i_t) \Delta W^i_t$
\STATE $b(t,\cdot)\!=\! \argmin_{b\in \Psi} \frac{1}{N}\sum_{i=1}^N \big[\frac{1}{2}\|b(X^i_{t})\|^2 \!+\! \trace(D(X^i_{t})\partial_x b(X^i_{t}))\big]$
\ENDFOR
\STATE $\cev X^i_T= X^i_T$
\STATE $\cev Y^i_T = g_f(\cev X^i_T)$, $\cev Z^i_T = \sigma(\cev X^i_T)^\top\partial_{x} g_f(\cev X^i_T)$,  and $\phi(T,\cdot)=g_f(\cdot)$
\FOR{$t \in\{T,T-\Delta t, \ldots,\Delta t\}$}
\STATE $\cev U^i_t = k(t,\cev X^i_t)$ and $\{\Delta \tilde W^i_t\}_{i=1}^N\sim N(0,\Delta t I_n )$
\STATE $\cev X^i_{t-\Delta t}=\cev X^i_t - a(\cev X^i_t,\cev U^i_t) \Delta t - b(t,\cev X^i_t)\Delta t -\sigma(\cev X^i_t)\Delta \tilde W^i_t$
\STATE $c(t,\cev X^i_t)=\trace(D(\cev X^i_t)\partial_{xx} \phi(t,\cev X^i_t)) -\partial_x \phi(t,\cev X^i_t)^\top b(t,\cev X^i_t)$
\STATE $\cev Y^i_{t-\Delta t}=\cev Y^i_t +g(t,\cev X^i_t,\cev Y^i_t,\cev Z^i_t) \Delta t + c(t,\cev X^i_t)\Delta t - (\cev Z^i_t)^\top \Delta \tilde W^i_t $
\STATE $\phi(t-\Delta t,\cdot)=\argmin_{\phi \in \Phi} \sum_{i=1}^N \| \cev Y^i_{t-\Delta t}-\phi(\cev X^i_{t-\Delta t})\|^2$
\STATE $\cev Z^i_{t-\Delta t} = \sigma(\cev X^i_{t-\Delta t})^\top \partial_x \phi(t-\Delta t,\cev X^i_{t-\Delta t})$

\ENDFOR
 \STATE \textbf{Output:}  $\{\phi(t,\cdot)\}_{t\in \{0,\Delta t,\ldots,T\} }$
\end{algorithmic}
\label{alg:TR}
\end{algorithm}

\section{NUMERICAL RESULTS}\label{sec:numerics}
We presented the LSMC and TR algorithms  for solving a general BSDE of the form~\eqref{eq:BSDE-general}. The output of these algorithms is used to compute the optimal control law for the SOC problem~\eqref{eq:SOC} according to the following iterative procedure. Starting from an initial  feedback control law $k(t,x)$ (typically chosen to be zero), the algorithms are implemented to obtain the function $\phi(t,x)$ in the relationship~\eqref{eq:Y-Z-phi}. The function $\phi(t,x)$ is then used to generate the next guess for the optimal control law, depending on which BSDE is solved. In particular,  
\begin{align*}
    \text{BSDE~\eqref{eq:BSDE-value-function}:}&~      k(t,x) = 
     - R^{-1} \tilde B(x)^\top \sigma(x)^\top \partial_x \phi(t,x),
    \\
    \text{BSDE~\eqref{eq:MP-BSDE}:}&~
     k(t,x)=\argmin_u H(x,u,\phi(t,x),\sigma(x)^\top\partial_x \phi(t,x)).
\end{align*}
Then, the procedure is iterated until the control law converges. 

Given the two BSDEs, \eqref{eq:BSDE-value-function} and \eqref{eq:MP-BSDE}, and two different approaches, LSMC and TR, to solve it, we end up with the following four methods to solve the SOC problem: 
\begin{itemize}
    \item LS-V: LSMC method for the value function BSDE~\eqref{eq:BSDE-value-function}
    \item LS-C: LSMC method for the co-state BSDE~\eqref{eq:MP-BSDE}
    \item TR-V: TR method for the value function BSDE~\eqref{eq:BSDE-value-function}
    \item TR-C: TR method for the co-state BSDE~\eqref{eq:MP-BSDE}
\end{itemize}

The algorithms are tested in the LQ setting with the following parametric class of functions, depending on which BSDE is solved: 
\begin{align*}
    \text{BSDE~\eqref{eq:BSDE-value-function}:}&\quad      \Phi =\{x \mapsto  
     \frac{1}{2}x^\top G x + g;\,(G,g) \in S^n\times \Re\}
    \\
    \text{BSDE~\eqref{eq:MP-BSDE}:}&\quad  
    \Phi =\{x \mapsto  
   G x;\,G \in \Re^{n\times n}\}
\end{align*}
where $S^n$ denotes the space of $n\times n$ symmetric positive definite matrices.  In the LQ setting, the probability distribution of the state $X_t$ is Gaussian. Consequently, the score function  admits the explicit formula 
\begin{align*}
    b(t,x) = D \Sigma_t^{-1}(X_t-m_t).
\end{align*}
where $m_t$ and  $\Sigma_t$ are the mean and covariance matrix for $X_t$. 
The affine form of the score function suggests to choose the affine function class \[\Psi=\{x\mapsto Ax+b;\,(A,b)\in \Re^{n\times n} \times \Re^n\}\] for the score function approximation~\eqref{eq:b-min} in the TR algorithm. 
In fact, in this case, the empirical minimization problem~\eqref{eq:b-min} admits an explicit solution $b(t,x) = D (\Sigma^{(N)}_t)^{-1}(X_t-m^{(N)}_t)$ where $m^{(N)}_t$ and $\Sigma^{(N)}_t$ are the empirical mean and covariance matrix formed from the samples $\{X^i_t\}_{i=1}^N$. This formula is used in our numerical implementation of the TR algorithm. The Python code used to generate the results appears in the  Github repository~\footnote{\url{https://github.com/YuhangMeiUW/LS-TR-BSDE}}. 

\subsection{Convergence and accuracy}\label{sec:numerics-accuracy}
The four algorithms are implemented on a two-dimensional LQ example 
with the model parameters  
\begin{align*}
&A = \begin{bmatrix}
	0 & 1 \\
	-1 & -0.1
\end{bmatrix},~B = \begin{bmatrix}
0\\1
\end{bmatrix},~\sigma = \begin{bmatrix}
1 & 0 \\
0 & 1
\end{bmatrix}
,\\&R=1,~m_0=\begin{bmatrix}
1\\0
\end{bmatrix},\quad Q=Q_f=\Sigma_0=\begin{bmatrix}
1 & 0 \\
0 & 1
\end{bmatrix}
\end{align*}
where $m_0$ and $\Sigma_0$ are the mean and covariance of the Gaussian initial distribution $p_0$. The time horizon $T=4$. The number of samples $N=2000$ and the time-discretization step-size $\Delta t = 0.02$. All four algorithms start with a zero control law $k(t,x)=0$. Each run of the algorithm results in a $2\times 2$ time-varying matrix $G_t$, which is used to update the control law according to the formula 
\begin{align*}
     k(t,x)=-R^{-1} B^\top G_t x. 
\end{align*}
The new control law is used to run the algorithm again, and this procedure is repeated 200 times to ensure convergence and fair comparison among all algorithms. 

Figure \ref{fig:g_matrices}  presents the four entries of the the resulting matrix $G_t$, after 200 iterations, in comparison to their exact values $G^*_t$ derived from the Riccati equation~\eqref{eq:Ricatti}.  The result shows that, although all four algorithms approximate the exact value reasonably, the TR-C algorithm admits much better accuracy compared to the other methods. In order to quantify this, 
we introduce the following mean-squared-error (MSE) criteria 
\begin{align}\label{eq:MSE} 
    MSE = \frac{1}{Tn^2}\int_0^T \|G_t^{output} - G_t^*\|_F^2 \ud t.
\end{align}
where $\|\cdot \|_F$ denotes the Frobenius norm of matrices, and we introduced the factor $\frac{1}{Tn^2}$ to normalize with respect to the time-horizon and the number of entries  in the matrix. The result of approximating the MSE, by averaging over 15 number of simulations and approximating the integral over the discrete-time horizon, appears in  Table~\ref{table:mse}, highlighting the significant accuracy of the TR-C method.  

\begin{table}[h]
\centering
\begin{tabular}{|c|c|}
\hline
Algorithm & MSE \\ \hline
LS-V      & $(4.5 \pm 1.9) \times 10^{-3}$    \\ \hline
LS-C      & $(4.8 \pm 2.5) \times 10^{-3}$   \\ \hline
TR-V      & $(6.1 \pm 1.5) \times 10^{-4}$    \\ \hline
TR-C      & $\mathbf{(2.2 \pm 0.4) \times 10^{-6}}$    \\ \hline
\end{tabular}
\caption{Numerical results for Section~\ref{sec:numerics-accuracy}. The MSE~\eqref{eq:MSE} evaluated for all four algorithms, averaged over 15 experiments.}
\label{table:mse}
\end{table}

\begin{figure}[thpb]
    \centering
    \includegraphics[width=0.6\linewidth]{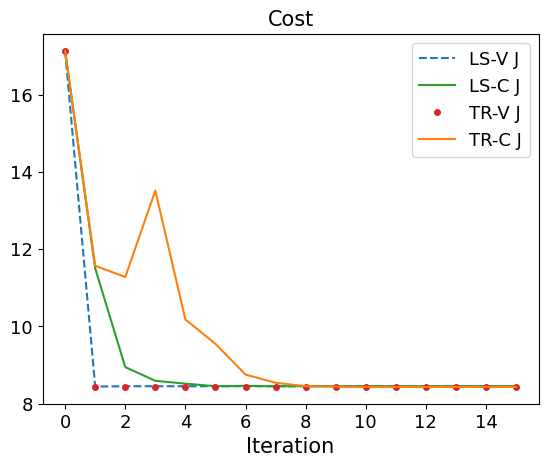}
    \caption{Numerical result for Section~\ref{sec:numerics-accuracy}: The value of the SOC cost~\eqref{eq:SOC} in the first 15 iterations of the four methods.}
    \label{fig:cost}
\end{figure}

The convergence of the SOC  cost function through the iterative implementation of the four algorithms is depicted in Figure \ref{fig:cost}. It is observed that the  LS-V and TR-V  algorithms both almost converge in one  iteration, while the  LS-C and TR-C algorithms take longer. This is due to the fact that,  solving the value-function BSDE~\eqref{eq:BSDE-value-function} results in the optimal control law~\eqref{eq:optimal-control-value-affine} irrespective of the control input that is used to generate the state. However,  solving the co-state BSDE~\eqref{eq:MP-BSDE} and updating the control input through~\eqref{eq:min-U} results  in the optimal control law when the control input that is used to generate the state is itself optimal.  

\subsection{Effect of time step-size and number of samples}\label{sec:dt-sample}
We study the influence of time step size $\Delta t$ and the number of samples $N$ on the accuracy of the four methods. In order to do so, we evaluate the MSE~\eqref{eq:MSE} averaged over 15 experiments for $\Delta t \in \{0.004, 0.02, 0.05, 0.1, 0.2, 0.3, 0.4\}$ and $N \in \{10, 50, 100, 500, 1000, 2000, 4000\}$. 
The modelling parameters is identical to the ones selected in Section~\ref{sec:numerics-accuracy}.


\begin{figure}[t]
    \centering
    \begin{subfigure}{0.235\textwidth}
        \centering
        \includegraphics[width=1\textwidth,trim={8pt 0pt 0pt 0 pt},clip]{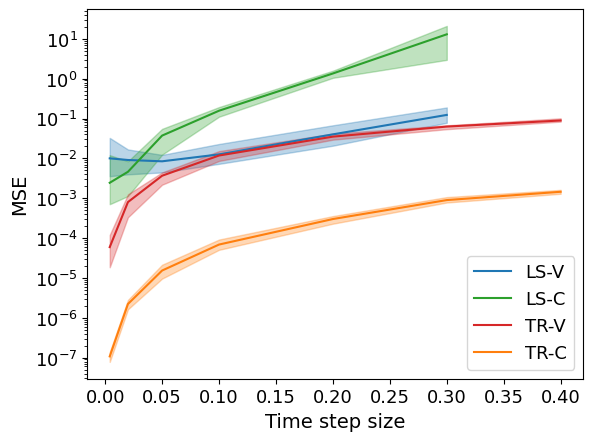}
        \caption{}
        \label{fig:mse_dt}
    \end{subfigure}
    \hfill
    \begin{subfigure}{0.235\textwidth}
        \centering
        \includegraphics[width=1\textwidth,trim={8pt 0pt 0pt 0 pt},clip]{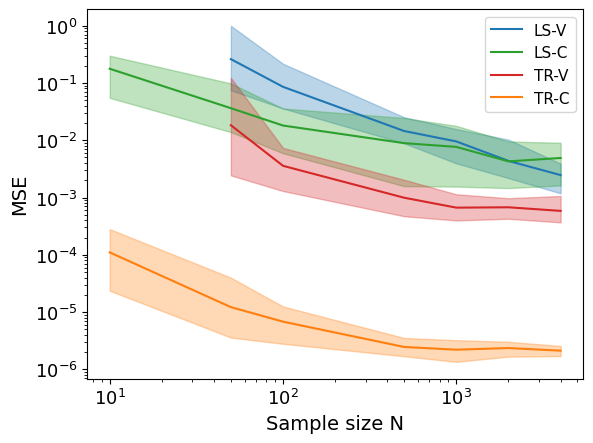}
        \caption{}
        \label{fig:mse_sample}
    \end{subfigure}
    \caption{Numerical results Sec. \ref{sec:dt-sample}. (a) MSE of four methods with $N=1000$ and different time step-size $\Delta t$;  The experiment is done with time horizon $T=4$, and sample size $N=1000$. (b) MSE of four methods with $\Delta t=0.02$ and different number of samples $N$. The shaded region represents the range from the minimum to the maximum across 15 experiments.}
\end{figure}

The result of the experiment with varying time step-size is presented in Figure~\ref{fig:mse_dt}.
It is observed that the MSE of all four methods increases with the time step size, as expected. However, it should be noted that, for the time step $\Delta t= 0.4$, the LS-V and LS-C methods experienced numerical instability 7 and 15 times, respectively, out of total 15 experiments.  The result shows that the LSMC approach is more sensitive to the step-size. 

The result of the experiment with varying number of samples is presented in Figure-\ref{fig:mse_sample}. 
It is also observed that the MSE of the four methods decreases as sample size $N$ increases. And, the LS-V and TR-V algorithm  experience numerical instability 11 and 15 times, respectively,  out of 15 experiments with a sample size of $10$. The result shows the sensitivity of solving value-function BSDE in comparison to co-state BSDE. 

The results from both experiments confirm the superior accuracy of the TR-C method, across different time step-size and number of samples. 

\subsection{Scalability with dimension}\label{sec:numerics-dim}
\begin{figure}[t]
    \centering
    \includegraphics[width=0.6\linewidth]{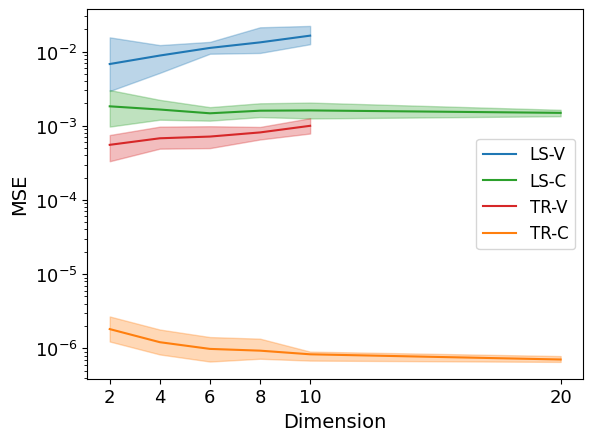}
    \caption{Numerical results for Section~\ref{sec:numerics-dim}: MSE of four methods as the problem dimension $n$ varies. Note that the MSE is normalized by the factor $n^2$. The shaded region represents the range from the minimum to the maximum across 15 experiments.}
    \label{fig:mse_dim}
\end{figure}

Our final experiment is concerned with scalability with the problem dimension. In order to do so, we use a mass-spring model with the system parameters
\begin{align*}
    A = \begin{bmatrix}
        0 & I_{p\times p} \\
        -\mathbb T & -I_{p\times p}
    \end{bmatrix}, ~ 
    B = \begin{bmatrix}
        0\\I_{p\times p}
    \end{bmatrix},\\ \Sigma = Q = Q_f = \Sigma_0 = I_{2p\times 2p},
    ~ R = I_{p\times p}.
\end{align*}
where  $p$ is the number of mass-springs in the system, $\mathbb T \in \R^{p \times p}$ is a Toeplitz matrix with 2 on the main diagonal and -1 on the first super and sub-diagonal, and $I_{p\times p}$ is $p\times p$ identity matrix. We performed 15 experiments with  time step-size $\Delta t=0.02$ and fixed sample size $N=1000$, as the dimension $n=2p$ varies. The result is depicted in Figure \ref{fig:mse_dim}.  It is observed that the MSE of two TR-C and LS-C methods, that are based on the co-state BSDE, to scale well with the problem dimension. However, the MSE of LS-V and TR-V algorithms  increases with the dimension. In fact both approaches experience numerical instability for dimension $n=20$. 

\section{CONCLUSIONS}
The paper presents a detailed exploration of two different numerical approaches, LSMC and TR, for solving the value function and co-state BSDEs that  arise in the context of stochastic optimal control. Through extensive numerical experiments, it was demonstrated that both methods are capable of solving the value function BSDE and the co-state BSDE, with the TR approach outperforming the LSMC method in terms of accuracy and numerical stability, particularly for the co-state BSDE. The results show that the TR method provides superior accuracy, specially in high-dimensional settings and under coarse time discretizations. A possible explanation for this observation is that the denoising procedure in the TR method reduces the noise in the regression problem for solving $\phi$ compared to the LSMC method. Theoretical justification for this observation is subject of future work. The efficiency of the TR approach, along with its robustness, suggests its potential applicability in a broader range of stochastic control problems beyond the linear-quadratic case, which is the subject of ongoing research.





 





\bibliographystyle{IEEEtran}
\bibliography{references}

\end{document}